\documentclass[12pt,a4paper]{amsart}
\usepackage[all]{xy}
\usepackage[latin1]{inputenc}
\usepackage{amsmath,hyperref}

\newcommand{\C}{{\mathbb C} }
\newcommand{\R}{{\mathbb R} }
\newcommand{\cA}{{\mathcal A} }

\newcommand{\cE}{{\mathcal E} }
\newcommand{\cF}{{\mathcal F} }

\newcommand{\cL}{{\mathcal L} }
\newcommand{\cM}{{\mathcal M} }

\newcommand{\cO}{{\mathcal O} }

\newcommand{\wh}{\widehat}
\newcommand{\wt}{\widetilde}
\def\ol#1{{\overline{#1}}}
\newtheorem{theorem} {Theorem} [section]
\newtheorem{definition}[theorem] {Definition}
\newtheorem{lemma}[theorem]{Lemma}
\newtheorem{remark}[theorem]{Remark}
\newtheorem{proposition}[theorem]{Proposition}
\newtheorem{corollary}[theorem]{Corollary}

\def\pw{Pe\-ters\-son-Weil }
\def\ks{Ko\-dai\-ra-Spen\-cer }
\def\ka{K{\"a}h\-ler}
\def\he{Her\-mite-Ein\-stein}

\def\ii{\sqrt{-1}}

\def\C{\mathbb{C}}

\def\cinf{C^\infty}
\def\rk{{\mathrm{rk}}}
\def\tr{{\mathrm{\, tr\, }}}
\def\ch{{\mathrm{ch}}}

\def\ab{{\alpha\ol\beta}}
\def\ba{{\ol\beta\alpha}}
\def\gba{{g^{\ol\beta\alpha}}}

\def\dbs{{\ol\partial^*}}
\def\db{{\ol\partial}}

\begin{document}

\title{Geometry of moduli spaces of Higgs bundles}

\author[I. Biswas]{Indranil Biswas}
\address{School of Mathematics, Tata Institute of Fundamental
Research, Homi Bhabha Road, Bombay 400005, India}

\email{indranil@math.tifr.res.in}

\author[G. Schumacher]{Georg Schumacher}

\address{Fachbereich Mathematik und Informatik,
Philipps-Universität Marburg, Lahnberge, Hans-Meerwein-Strasse,
D-35032 Marburg, Germany}

\email{schumac@mathematik.uni-marburg.de}


\date{}

\begin{abstract}

We construct a Petersson-Weil type K\"ahler form on the moduli
spaces of Higgs bundles over a compact K\"ahler manifold. A fiber
integral formula for this form is proved, from which it follows that
the Petersson-Weil form is the curvature of a certain determinant
line bundle, equipped with a Quillen metric, on the moduli space of
Higgs bundles over a projective manifold. The curvature of the
Petersson-Weil K\"ahler form is computed. We also show that, under
certain assumptions, a moduli space of Higgs bundles supports of
natural hyper-K\"ahler structure.

\end{abstract}

\maketitle

\section{Introduction}

A Higgs bundle over a compact Kähler manifold $X$ is a pair of the
form $(E,\varphi)$, where $E$ is a holomorphic vector bundle over
$X$ and $\varphi$ a holomorphic section of $End(E)\otimes
{\Omega}^1_X$ satisfying the integrability condition
$\varphi\wedge\varphi =0$. Higgs bundles over a compact Riemann
surface were introduced by Hitchin in \cite{hitchin} where he
constructed their moduli and investigated the global, as well as the
local, structures of the moduli space. One of the main results of
\cite{hitchin} was that a stable Higgs bundle admits a unique
Hermitian-Yang-Mills connection. Simpson, initiating the study of
Higgs bundles over compact Kähler manifolds of arbitrary dimension,
proved that a stable Higgs bundle admits a unique
Hermitian-Yang-Mills connection \cite{simpson0}. He also constructed
the moduli space of Higgs bundles over a complex projective manifold
\cite{simpson2}.

The aim in this article is to study the local geometry of a moduli
space of Higgs bundles from the point of view of the generalized \pw
geometry, which has been carried out for the moduli spaces of stable
vector bundles (cf.\ \cite{s-t, b-s2}). Here a moduli space is by
definition a reduced complex space. With additional effort, also
non-reduced moduli spaces can be investigated.

For any Higgs bundle $(E,\varphi)$ over a compact Kähler manifold
$X$ there is an associated complex of ${\mathcal O}_X$-modules
$$
D^\bullet\, :\,
D^0:= End(E)\rightarrow D^1 := End(E)\otimes \Omega^1_X \rightarrow
\cdots \rightarrow D^i := End(E)\otimes \Omega^i_X \rightarrow \cdots
$$
with the homomorphisms defined by $s\mapsto [s,\varphi]$.
The global endomorphisms and the infinitesimal deformations
of $(E,\varphi)$ are given by the hypercohomologies
${\mathbb H}^0(D^\bullet)$ and ${\mathbb H}^1(D^\bullet)$
respectively. Similarly the obstructions for deformations of
$(E,\varphi)$ are guided by
${\mathbb H}^2(D^\bullet)$. We consider the
Dolbeault resolution of the above complex $D^\bullet$.
So the hypercohomologies
of $D^\bullet$ get identified with the cohomologies associated
to the resulting double complex.

Now, the existence of Hermitian-Yang-Mills connections on stable
Higgs bundles allows us to introduce a natural inner product on the
terms of the above mentioned double complex, and general Hodge
theory provides the space of hypercohomologies with a hermitian
structure. In other words, it is possible to identify harmonic
representatives of the hypercohomologies, and the hermitian
structure on hypercohomologies is given by the hermitian structure
on the harmonic representatives. In particular, the space of
infinitesimal deformations of a stable Higgs bundle $(E,\varphi)$ is
equipped with a natural hermitian structure. In this way any moduli
space of Higgs bundles over $X$ is provided with a natural hermitian
metric.

This hermitian structure on a moduli space of Higgs bundles is
actually a Kähler structure (Proposition \ref{prop.ka}); we call
this Kähler form the generalized \pw form. We compute the curvature
tensor of this generalized \pw form (Theorem \ref{thm:curv}). For a
moduli space Higgs bundles over a compact Riemann surface, the
holomorphic sectional curvature turns out to be non-negative
(Corollary \ref{co1}).

Furthermore, we prove a fiber integral formula for the generalized
\pw form (see Proposition \ref{prop.fif} and Theorem \ref{th.fif}).
The fiber integral formula implies the Kähler property also for
families parameterized by singular base spaces. Finally, the
generalized Riemann-Roch theorem of Bismut, Gillet and Soul\'e
provides a certain determinant line bundle equipped with a Quillen
metric over the moduli space of Higgs bundles whose curvature form
coincides with the generalized \pw\ form on the moduli space
(Theorem \ref{thm3}).

In Section \ref{sec7},
we construct a distinguished locally exact holomorphic
$2$-form $\pi$ on the moduli space of Higgs bundles.
In order to show non-degeneracy of $\pi$,
we need an involution $\iota$, on the first hypercohomology, defined
in terms of harmonic representatives. To construct $\iota$ we need two
assumptions on $(E,\varphi)$: (1) the rational characteristic classes
of the projective bundle ${\mathbb P}(E)$ vanish, and (2)
$\dim {\mathbb H}^2(D^\bullet) = 1$.

We note that in general,
$\dim {\mathbb H}^2(D^\bullet) \geq 1$. Hence by semicontinuity,
the condition $\dim {\mathbb H}^2(D^\bullet) = 1$
defines a Zariski open subset of any
moduli space of Higgs bundles.
We also note that in general
$\dim {\mathbb H}^1(D^\bullet)$ may be odd.

As an application we construct a hyper-Kähler structure on the
moduli space using
$\iota$, the generalized \pw form and $\pi$ (Theorem \ref{tf}).

\section{Basic definitions}

Let $X$ be a compact, connected \ka\ manifold of dimension $n$
equipped with a \ka\ form $\omega_X$. We will write
$\omega_X = \ii g_{\alpha\ol \beta} dz^\alpha \wedge
dz^\ol\beta$ with respect to local
holomorphic coordinates $(z^1,\dots,z^n)$, and we will always
use the summation convention.

If $\cF$ is a coherent $\cO_X$-module, then the degree of $\cF$ with
respect to $\omega_X$ is defined as
$$
\deg\cF\, :=\, \int_X c_1(\cF)\wedge\omega^{n-1}\, .
$$
We denote by $E$ a holomorphic vector bundle over $X$ of rank $r$.
\begin{definition}
\begin{enumerate}
\item[(i)]
{\rm A {\it Higgs field} on a vector bundle $E$
over $X$ is a holomorphic section
$$
\varphi\in H^0(X,End(E)\otimes_{\cO_X}\Omega^1_X)
$$
such that
\begin{equation}\label{integrability}
\varphi\wedge\varphi =0,
\end{equation}
i.e., $[\varphi_\alpha,\varphi_\gamma]=0$ for all $\alpha, \gamma$,
where $\varphi=\sum_{\alpha=1}^n\varphi_\alpha dz^\alpha$.}
\item[(ii)]
{\rm A {\it Higgs bundle} is a pair $(E,\varphi)$, where $\varphi$
is a Higgs field on $E$.}
\end{enumerate}
\end{definition}

The definition of stability in this context is the following:

\begin{definition}
\begin{enumerate}
\item[(i)]
{\rm A Higgs bundle $(E,\varphi)$ is called {\it stable}, if
$$
\frac{\deg\cF}{\rk \cF}<\frac{\deg E}{\rk E}
$$
for all ${\mathcal O}_X$-coherent subsheaves $\cF$ of $E$ satisfying
the conditions \hfill \break $\varphi(\cF) \subset \cF
\otimes_{\cO_X}\Omega^1_X$ and $0< \rk \cF < \rk E$.}
\item[(ii)]
{\rm A {\it polystable} Higgs bundle is a direct sum of stable Higgs
bundles $(E^\nu,\varphi^\nu)$ with the same quotient $\frac{\deg
E^\nu}{\rk E^\nu}$.}
\end{enumerate}
\end{definition}

Polystable Higgs bundles $(E,\varphi)$
are known to carry a unique
Hermitian-Yang-Mills connection by results of Hitchin
and Simpson \cite{hitchin,simpson0}.

\begin{definition}\label{def.p.y}
{\rm Let $(E,\varphi)$ be a Higgs bundle. A
{\it Hermitian-Yang-Mills connection} on
$(E,\varphi)$ is a hermitian connection $\theta_E$ on $E$ with
curvature form $\Omega_E$ such that
\begin{equation}\label{hermein}
\Lambda(\Omega_E + \varphi\wedge\varphi^*)= \lambda \cdot {\rm id}_E
\end{equation}
for some $\lambda \in \R$, where $\Lambda=\Lambda_X$ is the adjoint
to the exterior multiplication of a form with $\omega_X$. In local
holomorphic coordinates $z^\alpha$ this equation reads
$$
g^{\ol\beta\alpha}\left(R_\ab+
[\varphi_\alpha,\varphi^*_\ol\beta]\right) = \lambda\cdot {\rm id}_E,
$$
where $\Omega_E = R_\ab\, dz^\alpha\wedge dz^\ol\beta$.}
\end{definition}

General theory provides a semi-universal deformation of pairs
$(E,\varphi)$, where $\varphi$ is a $End(E)$-valued holomorphic
$1$-form. The integrability condition $\varphi\wedge\varphi=0$
defines a complex analytic subspace of the parameter space, and thus
yields a semi-universal deformation for Higgs bundles. It follows
like in the classical case that {\em stable} Higgs bundles are {\em
simple} in the sense
\begin{equation}\label{eq.end}
H^0(X,End(E,\varphi))=\C\cdot {\rm id}_E\, ,
\end{equation}
where
\begin{equation}\label{eq.end2}
End(E,\varphi)\subset End(E)
\end{equation}
is the subsheaf that
commute with $\varphi$.
A {\em holomorphic family} $(E_s,\varphi_s)_{s \in S}$ of Higgs
bundles, parameterized by a complex space $S$, consists of a
holomorphic vector bundle $\cE$ on $X\times S$ and a
holomorphic section $\Phi$
of $End(\cE)\otimes \rho^*\Omega^1_X$, where $\rho: X\times S \to X$
is the canonical projection, such that $\cE| X\times \{s\} = E_s$,
and  $\Phi| X\times \{s\} =\varphi_s$ for all $s\in S$.

Observe that $\Phi$ defines an $End(E)$-valued, holomorphic $1$-form
on $X\times S$, as $\rho^*\Omega^1_X \subset \Omega^1_{X\times S}$.

Let $(E,\varphi)$ be any stable Higgs bundle.
In a local holomorphic family of Higgs bundles over a pointed space
$(S,s_0)$ with $(E,\varphi)$ as the central fiber (the Higgs
bundle over $s_0$), any isomorphism of the central fiber can be
extended to the restriction of the family
over a neighborhood of $s_0$. So stable Higgs bundles
possess universal deformations by general deformation theory
(cf.\ \cite{sch}).

As the uniquely determined
Hermitian-Yang-Mills connections on stable Higgs bundles
depend in a $\cinf$ way on the parameter of a holomorphic family of
such bundles, again by general results (even in the non-reduced
category) a {\em coarse moduli space} exists (cf.\ \cite{sch}).

We will denote by $\cM_H$ a moduli space of stable Higgs bundles
over $X$.

We will use the following conventions. The Kähler form $\omega_X$
gives rise to a connection on $X$, which we will extend in a flat
way to $X\times S$. As above, we will denote by $z^\alpha,
z^\gamma,\ldots$ local coordinates on $X$ together with the
conjugates $z^\ol\beta, z^\ol\delta,\ldots$, and by $s^i, s^k,
\ldots$ and $s^\ol\jmath, s^\ol\ell,\ldots$ respectively similar
coordinates on $S$. We use the semi-colon notation for covariant
derivatives of sections and differential forms or tensors with
values in the respective vector bundles induced by the Kähler metric
on $X$ and the hermitian connection on the bundle. Let the hermitian
connection $\theta_E$ on $E$ be given locally by $End(E)$-valued
$(1,0)$-forms $\{\theta_\alpha\}_{\alpha=1}^n$ with respect to some
local trivialization of $E$. Let $\sigma$ be a locally defined
section of $End(E)$, which is a matrix-valued function with respect
to the trivialization of $E$. We use
$$
\frac{\partial \sigma}{\partial z^\alpha}= \partial_\alpha \sigma =
\sigma_{|\alpha}
$$
and set
$$
\nabla_\alpha \sigma = \sigma_{;\alpha} = \sigma_{|\alpha} +
[\sigma, \theta_\alpha],
$$
and
$$
\sigma_{;\ol\beta}=\sigma_{|\ol\beta}.
$$
Hence
$$
\sigma_{;\ab} = \sigma_{;\ol\beta\alpha} + [\sigma, R_\ab],
$$
where $R_\ab$ denote the components of the curvature form
$\Omega_\ab= \theta_{\alpha|\ol\beta}$. For {\it tensors} with
values in the endomorphism bundle, we also have the contributions
that arise from the Kähler connection on the base. We denote by $g\,
dV$ the volume element $\omega_X^n/n!$  of the given \ka\ form.

\section{Infinitesimal deformations of Higgs bundles}\label{sec3}

Let $(E,\varphi)$ be any Higgs bundle over the compact Kähler
manifold $X$. For any integer $i\geq 0$, the Higgs field $\varphi$
gives a ${\mathcal O}_X$-linear homomorphism
$$
f_\varphi(i) \, :\, {End}(E)\otimes \Omega^i_X
\, \longrightarrow\, {End}(E)\otimes\Omega^{i+1}_X
$$
defined by $s\longmapsto [s,\varphi]$. From the given condition
that $\varphi\wedge\varphi\,=\, 0$ it follows immediately that
$$
f_\varphi(i+1)\circ f_\varphi(i) \, =\, 0
$$
for all $i$.
In other words, there is a complex of
${\mathcal O}_X$-coherent modules
$$
D^\bullet : 0 \to D^0:= End(E)
\stackrel{f_\varphi(0)}{\longrightarrow} End(E)\otimes \Omega^1_X
\stackrel{f_\varphi(1)}{\longrightarrow}\cdots
\stackrel{f_\varphi(i-1)}{\longrightarrow} D^i:= End(E)\otimes
\Omega^i_X
$$
$$
\stackrel{f_\varphi(i)}{\longrightarrow}
 End(E)\otimes\Omega^{i+1}_X
\stackrel{f_\varphi(i+1)}{\longrightarrow}\cdots
\stackrel{f_\varphi(n-1)}{\longrightarrow}
End(E)\otimes\Omega^{n}_X
\longrightarrow 0
$$
over $X$. We note that $H^0(X, End(E,\varphi))= {\mathbb
H}^0(D^\bullet)$ (see \eqref{eq.end2}). The space of all
infinitesimal deformations of $(E,\varphi)$ is parameterized by the
first hypercohomology ${\mathbb H}^1(D^\bullet)$, and the
obstructions to deformations of $(E,\varphi)$ are guided by
${\mathbb H}^2(D^\bullet)$; see \cite{biswas} for the details.

For computational convenience we
will work with the Dolbeault resolution of
the above complex $D^\bullet$.

Consider the spaces
$$
C^{p,q}\,:\, =\cA^{p,q}(X,End(E))
$$
of differentiable $(p,q)$-forms over $X$ with values in
$End(E)$ equipped with the Dolbeault operator
$$
d'': C^{p,q} \to C^{p,q+1}
$$
which is the $\ol\partial$-operator on $End(E)$-valued forms.
We also have an operator
$$
d' : C^{p,q} \to C^{p+1,q}
$$
which is defined by
$$
d'(\chi)=[\chi,\varphi]\, .
$$
Here the Lie bracket operation sends
$$
\chi_{\gamma\ol\delta} dz^\gamma\wedge dz^\ol\delta =
\chi_{\gamma_1,\dots,\gamma_p,\ol\delta_1,\dots,\ol\delta_q}
dz^{\gamma_1}\wedge\dots\wedge dz^{\gamma_p} \wedge
dz^{\ol\delta_1}\wedge\dots\wedge dz^{\ol\delta_q}
$$
to
$$
[\chi,\varphi]=[\chi_{\gamma\ol\delta},
\varphi_\alpha]dz^\gamma\wedge dz^\ol\delta\wedge dz^\alpha.
$$
Since the section $\varphi$ is holomorphic with $\varphi\wedge
\varphi =0$, it follows that $(C^{\bullet\bullet},
d',d'')$ is actually a double complex.
This double complex gives rise to a
degenerating spectral sequence, which converges to the
hypercohomology of the complex $D^\bullet$ defined earlier.

For the induced single
complex $(C^\bullet, d)$ with
\begin{equation}\label{de-cr}
C^r\,:=\, \bigoplus_{p+q=r}C^{p,q}
\end{equation}
we use the homomorphism
$$
d\,:=\, d'' + (-1)^{q+1} d'\, .
$$
The groups ${\mathbb H}^q(C^{\bullet\bullet}):= H^q(C^\bullet)$, for
$q=0,1$, are computed from the truncated complex
\begin{small}
\begin{equation}
\xymatrix{ 0 \ar[r] & C^{0,0}
\ar[r]^{\hspace{-2mm}{d^0}\hspace{3mm}} & C^{1,0}\oplus C^{0,1}
\ar[r]^{{\hspace{-2mm}d^1\hspace{3mm}}}  & C^{2,0} \oplus
C^{1,1}\oplus C^{0,2} }
\end{equation}
\end{small}
where
\begin{equation}\label{eb1}
d^0(f)= (-[f,\varphi],\ol\partial f)
\end{equation}
and
\begin{equation}\label{eb2}
d^1(a,b)=(-[a,\varphi],\ol\partial a +[b,\varphi],\ol\partial b)
\end{equation}
are defined above.

\begin{lemma}
Assume that $(E,\varphi)$ is equipped with a Hermitian-Yang-Mills
connection (see Definition \ref{def.p.y}). Then any holomorphic
section of $End(E,\varphi)$ (defined in \eqref{eq.end2}) is parallel
with respect to the induced connection.
\end{lemma}

\begin{proof}
Let $\sigma$ be a holomorphic section of  $End(E,\varphi)$. From
(\ref{hermein}),
$$
\gba([\sigma,R_\ab]+[\sigma,[\varphi_\alpha,\varphi^*_\ol\beta]])=0\, .
$$
Hence
\begin{eqnarray*}
\int_X \gba\tr \sigma_{;\alpha} \sigma^*_{;\ol\beta}  \,g\,dV
&=&- \int_X \gba\tr \sigma_{;\alpha\ol\beta} \sigma^*  \,g\,dV\\
&=&-\int_X \gba\tr [\sigma,R_\ab]\sigma^*  \,g\,dV
\end{eqnarray*}
Now $[\sigma,\varphi_\alpha]=0$ implies that the above integral equals
\begin{eqnarray*}
-\int_X
\gba\tr[\varphi_{\alpha},[\varphi^*_{\ol\beta},\sigma]]\sigma^*
\,g\,dV&=& -\int_X \gba\tr\sigma^*
[\varphi_{\alpha},[\varphi^*_{\ol\beta},\sigma]]\,g\,dV\\
& =&-\int_X \gba\tr[\sigma^*,\varphi_{\alpha}][\varphi_{\ol\beta}^*
,\sigma ] \,g\,dV\\ &\leq& 0.
\end{eqnarray*}
So the integral vanishes, and $\sigma_{;\alpha}=0$.
\end{proof}

{}From now on, we assume that the given Higgs bundle $(E,\varphi)$ is
{\it stable}. Therefore, it carries a hermitian metric satisfying
the Hermitian-Yang-Mills equation (see Definition \ref{def.p.y}).
This metric is unique up to a dilation by a globally constant scalar.

Obviously the space $H^0(C^\bullet)$ consists of those holomorphic
sections of $End(E)$ which commute with $\varphi$. The stability
condition of $(E,\varphi)$ implies
that any such section is a constant scalar multiple of the
identity automorphism of $E$ (see \eqref{eq.end}).

As mentioned earlier, the hypercohomology ${\mathbb H}^1(D^\bullet)=
\mathbb H^1(C^{\bullet\bullet})= H^1(C^\bullet)$ is the space of
all infinitesimal deformations of
$(E,\varphi)$. We denote the \ks map by
\begin{equation}\label{d-ks}
\rho\, :\, T_{s_0}S \to \mathbb H^1(C^{\bullet\bullet})\, .
\end{equation}
Now $(C^\bullet,d)$ becomes an {\em elliptic complex}, when equipped
with the inner products induced by the hermitian metric on $E$ and
the Kähler metric $\omega_X$ on $X$. In particular the formal adjoint
operators to $d^r$ are in fact adjoint.

More precisely, let $\sigma, \tau \in C^{0,0}$, then
$$
\langle \sigma,\tau\rangle = \int_X \tr (\sigma \tau^*) g \, dV\, ,
$$
where $\tau^*$ denotes the adjoint section with respect to the
hermitian metric on $E$. For $End(E)$-valued $(1,0)$-forms $\varphi
= \varphi_\alpha dz^\alpha$ and $\psi = \psi_\alpha dz^\alpha$ we
get
$$
\langle \varphi, \psi  \rangle = \int_X g^\ba \tr \, \varphi_\alpha
\psi^*_\ol\beta\, g\, dV.
$$
As usual $\Lambda$ denotes the operator that is adjoint to the
exterior multiplication of a form with $\omega_X$. We mention we
have $[a,b]^*=-[a^*,b^*]$   for any forms $a,b \in C^{p,q}$.

For the computation of adjoint derivatives we need the following
notation: For $v \in C^{1,1}$ we consider $[v,\varphi^*]$ as a
tensor rather than as an alternating form, then
$$
\wt\Lambda [v_{\alpha\ol\beta},\varphi^*_\ol\delta]:=-
g^{\ol\delta\alpha}[v_{\alpha\ol\beta},\varphi^*_\ol\delta]dz^\ol\beta.
$$
So the contraction $\widetilde\Lambda[v,\varphi^*]$ stands for the
contraction of $v$ and $\varphi^*$, and it does not comprise
$\Lambda v$.

\begin{lemma}\label{adjoint}
For $(a,b) \in C^{1,0}\oplus C^{0,1}$ and $(u,v,w)\in C^{2,0}\oplus
C^{1,1}\oplus C^{0,2}$ the following hold:
\begin{gather}
d^{0*}(a,b)= -\Lambda[a,\varphi^*]+ \ol\partial^*b \label{dstar0}\\
d^{1*}(u,v,w)= (\Lambda[u,\varphi^*] +\ol\partial^*v,
\widetilde\Lambda[v,\varphi^*]+ \ol\partial^*w). \label{dstar1}
\end{gather}
where $d^0$ and $d^1$ are defined in \eqref{eb1} and \eqref{eb2}
respectively.
\end{lemma}

\begin{proof}
The first equation follows from $\langle (a,b), d f\rangle = -
\langle a,[f,\varphi]\rangle + \langle \dbs b,f\rangle$ and
$$-
\langle a,[f,\varphi]\rangle = \langle a, [f,\varphi]\rangle = -
\langle \Lambda[a,\varphi^*],f\rangle\, ,
$$
whereas for (\ref{dstar1}) we have
$$
\langle (u,v,w), (-[a,\varphi], \db a + [b,\varphi],\db
b)\rangle
$$
$$=\, \langle -\Lambda [u, \varphi^*],a\rangle  + \langle
\dbs v, a\rangle + \langle -\wt\Lambda[v,\varphi^*],b\rangle +
\langle \dbs v, a\rangle
$$
finishing the proof of the lemma.
\end{proof}

Let $(\cE,\Phi)$ be a holomorphic family of Higgs bundles over a
complex space $S$, and denote by $\{h_s\}$ any $\cinf$ family of
hermitian metrics on $\cE_s$, i.e.,\ a hermitian metric $h$ on $\cE$
over $X\times S$. Let $(s^1,\dots,s^k)$ be holomorphic coordinates
on $S$, if $S$ is smooth, or holomorphic coordinates on an ambient
smooth space into which a neighborhood of $s_0\in S$ is minimally
embedded.

Let $\Omega$ be the curvature form of the hermitian connection
for $h$ on $\cE$ over $X \times S$. The curvature tensor for this
connection will be denoted by $R$. So the contraction
$$
\Omega\; \llcorner \frac{\partial}{\partial s_i}
$$
equals
$$
R_{i\ol\beta}\, dz^\ol\beta\, .
$$

Following the construction in \cite{s-t}
one can see that the global tensors $\Omega$ and $\varphi$ over $X
\times S$ already describe the infinitesimal deformations.
In other words, we have the following lemma:

\begin{lemma}\label{lab}
The \ks class
$$
\rho\left(\left.\frac{\partial}{\partial s_i}\right|_{s_0}  \right)
\in \mathbb H^1(C^{\bullet\bullet})
$$
(the homomorphism $\rho$ is defined in \eqref{d-ks}) is represented by
\begin{equation}\label{defetai}
\eta_i\,=\,(\varphi_{\alpha;i}\,dz^\alpha\, ,
R_{i\ol\beta}\,dz^\ol\beta)|_{X\times\{s_0\}}.
\end{equation}
\end{lemma}

The \ka\ form $\omega_X$ and the hermitian metric $h$ together
provide the above double complex $C^{\bullet\bullet}$ with a
natural inner product such that the adjoint operators $d^{j*}$ are
the formal adjoint operators.

Now assume that for each point $s\in S$, the Higgs bundle
$(\cE_s,\Phi_s)$ over $X$ is stable. The Hermitian-Yang-Mills
connections on this family of stable Higgs  bundles $(\cE,\Phi)$
are induced by a hermitian metric $h$ on $\cE$, whose curvature form
$\Omega$ is unique up to a differential form of type ${\rm id}_E\otimes
f^*\omega'$, where $\omega'$ is some $(1,1)$-form on the base $S$
and $f : X\times S\rightarrow S$ is the natural projection. Indeed,
this follows immediately from the fact that any two
Hermitian-Yang-Mills metrics on a stable Higgs bundle differ by
multiplication with a constant scalar.

Therefore, the components $R_{i\ol\beta}$ of the curvature tensor in
Lemma \ref{lab} are uniquely determined by the family of
Hermitian-Yang-Mills connections on the Higgs bundles
$(\cE_s,\varphi_s)$.

\begin{proposition}\label{etaharm}
The $End(\cE)$-valued $1$-forms
$$
\eta_i=\varphi_{\alpha;i}dz^\alpha + R_{i\ol\beta}\, dz^\ol\beta
$$
are the harmonic representatives of the \ks\ classes
$\rho(\partial/\partial s^i|_{s_0})$.
\end{proposition}
\begin{proof}
Following Lemma  \ref{adjoint} we find
$$
d^*(\varphi_{\alpha;i}dz^\alpha,R_{i\ol\beta}\, dz^\ol\beta)= -
\Lambda[\varphi_{;i},\varphi^*] + \dbs(R_{i\ol \beta}dz^{\ol\beta})
= g^\ba(-[\varphi_{\alpha;i},\varphi^*_\ol\beta] - R_{\ab;i}) = 0
$$
because of the Hermitian-Yang-Mills condition (\ref{hermein})
for Higgs bundles.
\end{proof}

For applications in Section~\ref{nonabhodge} we introduce a
decomposition of the complex $D^\bullet$.

Let
\[
ad(E)\subset End(E)
\]
be the subbundle of trace zero
endomorphisms, and let
\[
pr : End(E) \to ad(E)
\]
defined by $pr(\chi):=
\chi - (1/\text{rk}(E)) \text{tr}(\chi) \text{id}_E$
be the projection onto
the trace free part. We extend the homomorphism $pr$
to the complex $D^\bullet$ and its resolution. Now the
complex of ${\mathcal O}_X$-coherent modules
$$
D^\bullet_0 : 0 \to D^0_0:= ad(E)
\stackrel{f_\varphi(0)}{\longrightarrow} ad(E)\otimes \Omega^1_X
\stackrel{f_\varphi(1)}{\longrightarrow}\cdots
\stackrel{f_\varphi(i-1)}{\longrightarrow} D^i_0:= ad(E)\otimes
\Omega^i_X
$$
$$
\stackrel{f_\varphi(i)}{\longrightarrow}
 ad(E)\otimes\Omega^{i+1}_X
\stackrel{f_\varphi(i+1)}{\longrightarrow}\cdots
\stackrel{f_\varphi(n-1)}{\longrightarrow} ad(E)\otimes\Omega^{n}_X
\longrightarrow 0
$$
over $X$, is a sub-complex of the complex $D^\bullet$.

We identify
$\Omega_X^k$ with $\text{id}_E\otimes\Omega_X^i \subset D^i$.

\begin{lemma}\label{decomp}
The restrictions of the chain morphisms $d^k$ to $\Omega^k_X$ are
identically zero. Moreover,
$$
D^\bullet= D^\bullet_0 \oplus \Omega^\bullet_X
$$
is an {\it orthogonal} decomposition. The resolution of
$C^{\bullet\bullet}$ also decomposes in an orthogonal way
$$
C^{\bullet\bullet}=C^{\bullet\bullet}_0 \oplus \mathcal
A_X^{\bullet\bullet},
$$
where $\mathcal A_X^{\bullet\bullet}$ corresponds to the Dolbeault
resolution of $\Omega_X^\bullet$. In particular the following holds.

The sub-complex $C^{\bullet\bullet}_0 \subset C^{\bullet\bullet}$
is preserved by both $d_0$ and $d_0^*$.
Let $d_0$ and $d_0^*$ be the restrictions of $d$ and $d^*$
respectively to $C^{\bullet\bullet}_0$. Then
$$
pr \circ d = d_0 \circ pr \text{  and } pr \circ d^* = d_0^* \circ
pr.
$$
\end{lemma}

\begin{proof}
The proof follows from $\tr[\chi,\varphi]=0$ for any $\chi \in D^k$,
and the simple fact that covariant derivatives commute with taking
traces.
\end{proof}

Concerning the second cohomology, we note
\begin{lemma}\label{lale}
There is a natural embedding
$$
\mathbb C\cdot\omega_X\cdot {\rm id}_E \hookrightarrow {\mathbb
H}^2(C^{\bullet\bullet})\, .
$$
\end{lemma}

\begin{proof}
We have to consider
$$
\epsilon=(0, \omega_X \cdot {\rm id}_E ,0) \in C^{2,0}\oplus C^{1,1}\oplus
C^{0,2}\, .
$$
Then $d\epsilon=(0,[\omega_X\cdot {\rm id}_E,\varphi],\db (\omega_X \cdot
{\rm id}_E),0)=0$, and from \eqref{dstar1},
$$
d^*\epsilon = (\db^*(\omega_X\cdot {\rm id}_E) \, , \wt\Lambda[\omega_X
\cdot {\rm id}_E,\varphi^* ])=0\, .
$$
So $\epsilon$ is harmonic, and it is different from zero, since $\int_X
\tr\epsilon^n\neq 0$.
\end{proof}

We will need the following proposition.

\begin{proposition}
The obstructions for deformations of a Higgs bundle $(E,\varphi)$
are already contained in the second hypercohomology
${\mathbb H}^2(D^\bullet_0)$. If $\dim {\mathbb H}^2(D^\bullet) =1$,
then we have ${\mathbb H}^2(D^\bullet_0)=0$.
\end{proposition}

\begin{proof}
Note that the deformations of any Higgs line bundle are unobstructed
(as the deformations of a line bundle are so). Therefore, setting
$G=\text{PGL}(r, {\mathbb C})$ in Theorem 3.1 of \cite{br}, where
$r= \text{rank}(E)$, we conclude that if the image of
${\mathbb H}^2(D^\bullet)$ in ${\mathbb H}^2(D^\bullet_0)$ is zero,
then all deformations of $(E,\varphi)$ are unobstructed.

If $\dim {\mathbb H}^2(D^\bullet) =1$, then ${\mathbb H}^2(D^\bullet)$
must be generated by the image of the
nonzero homomorphism in Lemma \ref{lale},
hence ${\mathbb H}^2(D^\bullet_0)$ vanishes.
\end{proof}

\section{Generalized \pw\ metric}

Now, we are in a position to introduce a generalized \pw metric
on the parameter space $S$ for a family of stable Higgs bundles.
The generalized \pw metric is an
inner product $G^{PW}$ on the tangent spaces $T_sS$ of the bases of
holomorphic families, which is positive definite for effective
families, and it is defined in terms of the tensors $\eta_i$
representing the \ks\ classes. This is possible, also in
the case where $S$ is
singular, because the family of Higgs forms, and the curvature
form for the connection on vector bundles, still
exist on the first infinitesimal neighborhood. The latter fact
follows from the approach that is based on the implicit function
theorem.

We call this hermitian structure the {\em generalized \pw metric},
and use the following notation:
\begin{gather}
G^{PW}\left(\left. \frac{\partial}{\partial s^i}\right|_{s_0},
\left. \frac{\partial}{\partial s^\ol\jmath}\right|_{s_0} \right) :=
G_{i\ol\jmath}^{PW}:= \langle\eta_i,\eta_j\rangle = \\
\hspace{2cm}= \int_X \tr (g^\ba
\varphi_{\alpha;i}\varphi^*_{\ol\beta,\ol\jmath})g\, dV + \int_X
\tr( g^\ba R_{i\ol\beta} R_{\alpha\ol\jmath}) g\, dV\, . \nonumber
\end{gather}
We set
$$
\omega_{PW}= \ii \; G_{i\ol\jmath}^{PW}ds^i\wedge ds^\ol\jmath\, .
$$
In order to compute the induced connection, we need certain
identities.

\begin{lemma}\label{dstaretaik}
Let $\eta_i=(\varphi_{\alpha;i}dz^\alpha, R_{i\ol\beta}dz^\ol\beta)$. Then
\begin{gather}
\eta_{i;k}=\eta_{k;i}\label{etasymm}\\
d\eta_{i;k}  + \eta_i\wedge \eta_k= 0 \label{deta_ik}\\
d^*\eta_{i;k}=0\label{formetastarik}\\
\eta_{i;\ol\jmath} = d R_{i\ol\jmath}\label{dRij}\\
\Box R_{i\ol\jmath}= d^*d
R_{i\ol\jmath}=d^*\eta_{i;\ol\jmath}\label{boxR}
\\
d^*\eta_{i;\ol\jmath}= g^\ba
([\varphi_{\alpha;i},\varphi^*_{\ol\beta,\ol\jmath}]
+[R_{i\ol\beta},R_{\alpha\ol\jmath}] ) \label{dstaretaij}
\end{gather}
\end{lemma}
\begin{proof}The symmetry \eqref{etasymm} of
$\eta_{i;k}$ follows immediately from
\hfill\break
$$
\eta_{i;k}\,=\, (\varphi_{\alpha;ik}dz^\alpha\, ,
R_{i\ol\beta;k}dz^\ol\beta )\, ,
$$
which is symmetric in $i$ and $k$.

We show (\ref{deta_ik}):
\begin{gather*}
d\eta_{i;k}=\hspace{9cm}\\(-[\varphi_{\alpha;ik},
\varphi_\gamma]dz^\alpha\wedge dz^\gamma ,
\varphi_{\alpha;ik\ol\beta} dz^\ol\beta\wedge dz^\alpha +
[R_{i\ol\beta;k} dz^\ol\beta,\varphi_\alpha dz^\alpha] ,
R_{i\ol\beta; k\ol\delta}dz^\ol\delta\wedge dz^\ol\beta )\\ =
([\varphi_{;i}, \varphi_{;k}], - [\varphi_{;i},
R_{k\ol\beta}dz^\ol\beta]- [\varphi_{;k}, R_{i\ol\beta}dz^\ol\beta],
- [R_{i\ol\beta}dz^\ol\beta, R_{k\ol\delta} dz^\ol\delta]).
\end{gather*}
This gives the formula (with the given sign convention for
$C^{\bullet\bullet}$).

Concerning (\ref{formetastarik}), we have
$$
d^*\eta_{i;k} = -g^\ba
[\varphi_{\alpha;ik},\varphi^*_\ol\beta]-g^\ba
R_{i\ol\beta;k\alpha}.$$
We use (\ref{hermein}) in \hfill\break
\begin{gather*}
-g^\ba R_{i\ol\beta;k\alpha}= -g^\ba R_{\ab;ik}=
g^\ba[\varphi_\alpha, \varphi^*_\ol\beta]_{;ik}\\= g^\ba
[\varphi_{\alpha;i},\varphi^*_\ol\beta]_{;k}=g^\ba
[\varphi_{\alpha;ik},\varphi^*_\ol\beta].
\end{gather*}

Next we show (\ref{dRij}),
\begin{gather*}
\eta_{i;\ol\jmath}= (\varphi_{\alpha;i} dz^\alpha\, ,
R_{i\ol\beta}dz^\ol\beta)_{;\ol\jmath}\\=((\varphi_{\alpha;\ol\jmath
i} +[\varphi_\alpha,R_{i\ol\jmath}])dz^\alpha\, ,
R_{i\ol\jmath;\ol\beta}dz^\ol\beta)\\=(-[R_{i\ol\jmath},\varphi]\, ,
\db R_{i\jmath})= dR_{i\ol\jmath}\, .
\end{gather*}

Formula (\ref{boxR}) is immediate. For the last formula
(\ref{dstaretaij}), we consider
\begin{gather*}
d^*\eta_{i;\ol\jmath} = g^\ba
([\varphi_{\alpha;i\ol\jmath},\varphi^*_\ol\beta] -
R_{i\ol\beta;\ol\jmath\alpha})=\hspace{4cm}\\
-g^\ba([\varphi_{\alpha;\ol\jmath i} + [\varphi_\alpha,
R_{i\ol\jmath}],\varphi^*_\ol\beta]+
(R_{i\ol\beta;\alpha\ol\jmath}-[R_{i\ol\beta},R_{\alpha\ol\jmath}]
))\, .\\
\end{gather*}
Using the Hermitian-Yang-Mills equation, we get
$$
g^\ba(-[ [\varphi_\alpha, R_{i\ol\jmath}],\varphi^*_\ol\beta]+
[\varphi_{\alpha;i\ol\jmath},\varphi^*_\ol\beta]+
[\varphi_{\alpha;i},\varphi^*_{\ol\beta;\ol\jmath}]+
[R_{i\ol\beta},R_{\alpha\ol\jmath}])=
g^\ba([\varphi_{\alpha;i},\varphi^*_{\ol\beta;\ol\jmath}]+
[R_{i\ol\beta},R_{\alpha\ol\jmath}])
$$
completing the proof of the lemma.
\end{proof}
As a consequence of Proposition~\ref{etaharm} and (\ref{dRij}), we
note that
\begin{equation}\label{orth_i_jk}
\langle\eta_i, \eta_{j;\ol k}\rangle= \langle \eta_i, dR_{j\ol
k}\rangle = \langle d^*\eta_i, R_{j\ol k}\rangle=0
\end{equation}
Using the above notation, we set
$$
G_{i\ol\jmath|k}^{PW}= \partial_{s_k}G_{i\ol\jmath}^{PW}\, .
$$

\begin{proposition}\label{prop.ka}
The generalized \pw metric is \ka, more precisely, we have
\begin{equation}\label{Gijk}
G_{i\ol\jmath|k}^{PW}\,=\,\langle \eta_{i;k}\, , \eta_j\rangle.
\end{equation}
\end{proposition}
\begin{proof}
We use $G_{i\ol\jmath|k}^{PW}\,=\,\langle \eta_{i;k}\, ,
\eta_j\rangle + \langle\eta_i, \eta_{j;\ol k}\rangle$, and
\eqref{orth_i_jk}.
\end{proof}

\begin{corollary}\label{cor:normalcoord}
Let $s_0$ be some point of the base $S$ of a universal deformation
of a Higgs-bundle. Consider normal coordinates $\{s^i\}$ for the \pw
metric at the base point
$s_0$. Then for all $i$, $k$, the harmonic projections
$H(\eta_{i;k}|_{s=s_0})$ vanish.
\end{corollary}
\begin{proof}
This follows immediately from the fact that the $\eta_i$ are
harmonic and span the whole space $H^1(C^\bullet)$.
\end{proof}

\section{Curvature of the generalized \pw metric}

Let $S=\{(s^1,\ldots,s^N)\}$ be the smooth base of a universal
deformation of a Higgs bundle equipped with a family of
Hermitian-Yang-Mills metrics.  Let $\eta_i=
(\varphi_{\alpha;i}dz^\alpha, R_{i\ol\beta}dz^\beta)$ be the
harmonic representative of the \ks class $\rho(\partial/\partial
s_i|_s) \in \mathbb H^1(C^{\bullet\bullet})$. We consider the
associated single complex $C^\bullet$ as an elliptic complex
equipped with the Laplacians $\Box= d^*d+ d d^*$ acting on
$End(E)$-valued forms in all degrees. This elliptic complex
possesses harmonic projections $H$ and Green's operators $G$.

\begin{theorem}\label{thm:curv}
Let
$\eta_i=(\varphi_{\alpha;i}dz^\alpha, R_{i\ol\beta}dz^\ol\beta)$
be the elements of a basis of
the harmonic \ks forms, depending on $s \in S$. Then,
 the curvature tensor of the generalized \pw metric equals
\begin{eqnarray}\label{eq111}
R^{PW}_{i\ol\jmath k\ol\ell} & = & \int_X \tr\left( G(
\Lambda(\eta_i\wedge\eta^*_\ol\jmath))
\Lambda(\eta_k\wedge\eta^*_\ol\ell)\right) \;g\,dV\label{eq:curv}\\
&& +\int_X \tr \left(G( \Lambda(\eta_i\wedge\eta^*_\ol\ell))
\Lambda(\eta_k\wedge\eta^*_\ol\jmath)\right) \;g\,dV\nonumber\\
&& + \int_X \tr\;\left([\eta_i \wedge \eta_k]\; \wedge
G([\eta_\ol\jmath^* \wedge \eta^*_\ol\ell])\right)
\frac{\omega_X^{n-1}}{(n-1)!}\, . \nonumber
\end{eqnarray}
Explicitly, we have
\begin{eqnarray}\label{eq112}
R^{PW}_{i\ol\jmath k\ol\ell} & = & + \int_X \tr(R_{i\ol\jmath}\Box
R_{k\ol\ell}+ R_{i\ol\ell}\Box R_{k\ol\jmath} )\;g\,dV\\
&& + \int_X \tr\left([\eta_i \wedge \eta_k]\; \wedge
G([\eta_\ol\jmath^* \wedge \eta^*_\ol\ell])\right)
\frac{\omega_X^{n-1}}{(n-1)!}\, . \nonumber
\end{eqnarray}
Here $[\eta_i\wedge\eta_k]$ and $[\eta_i\wedge\eta_\ol\jmath]$ stand
for the exterior product of forms with values in an endomorphism
bundle combined with the Lie product.
\end{theorem}
\begin{remark}
The first two terms in \eqref{eq111} and the first term in
\eqref{eq112} resp.\ are semi-positive. Because of the different
order of non-conjugate and conjugate terms in the first and second
part of the curvature formula, the last terms of \eqref{eq111} and
\eqref{eq112} resp.\ yield semi-negativity.
\end{remark}

\begin{proof}[Proof of Theorem \ref{thm:curv}]
We use normal coordinates at a given point of $S$. Then
\begin{equation}\label{d.A.B.}
-R^{PW}_{i\ol\jmath k\ol\ell} = G^{PW}_{i\ol\jmath|k
\ol\ell}=\langle \eta_{i;k\ol\ell}, \eta_j\rangle + \langle
\eta_{i;k},\eta_{j;\ell}\rangle=: A + B\, .
\end{equation}
We compute $A$.
$$
\eta_{i;k\ol\ell}=
(\varphi_{\alpha;ik\ol\ell}dz^\alpha,R_{i\ol\beta;k\ol\ell}dz^\ol\beta)\in
C^{1,1}.
$$
First we need
$$
\varphi_{\alpha;ik\ol\ell}= \varphi_{\alpha;i\ol\ell k} +
[\varphi_{\alpha;i},R_{k\ol\ell}]=
$$
$$
\varphi_{\alpha;\ol\ell ik}+[\varphi_\alpha,
R_{i\ol\ell}]_{;k}+[\varphi_{\alpha;i},R_{k\ol\ell}]=
[\varphi_{\alpha;k},R_{i\ol\ell}]+
[\varphi_{\alpha;i},R_{k\ol\ell}]+[\varphi_{\alpha},R_{i\ol\ell;k}]
$$
so that
$$
\langle\varphi_{\alpha;ik\ol\ell},\varphi_{\gamma;j}\rangle=\int_X
g^\ba \tr\left(([\varphi_{\alpha;k},R_{i\ol\ell}] +
[\varphi_{\alpha;i},R_{k\ol\ell}]+[\varphi_\alpha , R_{i\ol\ell;k} ]
) \varphi^*_{\ol\beta;\ol\jmath}\right)\; g\, dV
$$
and
$$
R_{i\ol\beta;k\ol\ell}= R_{i\ol\ell,k\ol\beta} -
[R_{i\ol\ell},R_{k\ol\beta}]+[R_{i\ol\beta} , R_{k\ol\ell}]
$$
so that
\begin{gather*}
\int_X  g^\ba\tr( R_{i\ol\ell;k\ol\beta} R_{\alpha\ol\jmath})\;g\,dV
= -\int_X g^\ba\tr(R_{i\ol\ell;k}R_{\ab;\ol\jmath})\;g\,dV
=  \hspace{2cm}\\
=
\int_Xg^\ba\tr(R_{i\ol\ell;k}[\varphi_{\alpha},
\varphi^*_{\ol\beta\ol\jmath}])\;g\,dV=
\int_Xg^\ba\tr([R_{i\ol\ell;k},
\varphi_{\alpha}]\varphi^*_{\ol\beta\ol\jmath})\;g\,dV,
\end{gather*}
which is inserted into the expression for
$$
\langle R_{i\ol\beta;k\ol\ell}dz^\ol\beta ,
R_{j\ol\delta}dz^\ol\delta\rangle.
$$
So far, we have
$$
A= \int_X
g^\ba\tr\left(R_{i\ol\ell}([R_{\alpha\ol\jmath},R_{k\ol\beta}]-
[\varphi_{\alpha;k},\varphi^*_{\ol\beta;\ol\jmath}]) +
R_{k\ol\ell}([R_{\alpha\ol\jmath},R_{i\ol\beta}]-
[\varphi_{\alpha;i},\varphi^*_{\ol\beta;\ol\jmath}]) \right)
\;g\,dV\, ,
$$
where $A$ is defined in \eqref{d.A.B.}.

Now, we compute the $(1,1)$-component of
$\eta_k\wedge\eta_\ol\jmath^*$:
$$
-(\eta_k\wedge\eta_\ol\jmath^*)_{(1,1)}= (
[R_{\alpha\ol\jmath},R_{k\ol\beta}]-[\varphi_{\alpha
;k},\varphi^*_{\ol\beta,\ol\jmath}] )dz^\alpha \wedge dz^\ol\beta
$$
so that with (\ref{boxR}) and (\ref{dstaretaij}),
\begin{equation}
\Lambda(\eta_i\wedge\eta^*_\ol\jmath)= d^*\eta_{i;\ol\jmath}= \Box
R_{i\ol\jmath}.
\end{equation}
Hence
\begin{eqnarray*}
A&=& - \int_X \tr(R_{i\ol\jmath}\Box R_{k\ol\ell}+ R_{i\ol\ell}\Box
R_{k\ol\jmath}  )\;g\,dV\\
&=& - \int_X \tr\left(G\left(\Lambda (\eta_i\wedge
\eta^*_\ol\jmath)\right) \Lambda(\eta_k \wedge \eta^*_\ol\ell )+
G\left(\Lambda (\eta_i\wedge \eta^*_\ol\ell)\right) \Lambda(\eta_k
\wedge \eta^*_\ol\jmath ) \right)\;g\,dV
\end{eqnarray*}

We compute $B$ defined in \eqref{d.A.B.}.
Since $d^*\eta_{i;k}=0$, and $H(\eta_{i;k})=0$ by
(\ref{formetastarik}) and Corollary~\ref{cor:normalcoord}, we have
$$
\eta_{i;k}= G d^*d  \eta_{i;k}= d^*G d \eta_{i;k},
$$
and
$$
B= \langle d\eta_{i;k}\, , G d \eta_{j;\ell}\rangle
$$
so that we get the third term of (\ref{eq:curv}) using
(\ref{deta_ik}). This completes the proof of the theorem.
\end{proof}

We estimate the holomorphic sectional curvature for $\dim X=1$:
$$
R_{i\ol i i\ol i}^{PW} = 2 \langle dR_{i\ol i}\, ,dR_{i\ol i}\rangle
\geq 0\, .
$$
with $d=d^0$ (see \eqref{eb1}). Equality holds only if
$dR_{i\ol\jmath}=0$, that is, $R_{i\ol i}$ is a holomorphic section of
$End(E,\varphi)$; see \eqref{eq.end2} for the
definition of $End(E,\varphi)$. Since
$(E,\varphi)$ is stable, any holomorphic section of
$End(E,\varphi)$ is a constant multiple of the identity
automorphism of $E$ (see \eqref{eq.end}).

Therefore, we have the following corollary:

\begin{corollary}\label{co1}
When $\dim X =1$, the holomorphic sectional curvature of the \pw\
metric is non-negative.
\end{corollary}

\section{Fiber integral formula}

We will show the existence of a local
$\partial\ol\partial$-potential for the generalized \pw\ metric on a
base space $S$ of a universal deformation. We note that this
implies the \ka\ condition of the \pw\ metric.

We consider a moduli space of stable Higgs bundles $\cM_H$.
Although, in general there is no universal holomorphic vector bundle
$\cE$ globally on $X\times\cM_H$, the bundle $End(\cE)$ exists in
the {\it orbifold sense} over all of $X \times \cM_H$, since the
non-zero scalars act trivially on $End(\cE)$. We will furthermore
need the highest exterior power $\Lambda^r\cE$, a tensor power of
which also descends to $X\times \cM_H$.

Representing a point $p \in \cM_H$ by the isomorphism class of a
Higgs bundle $(E,\varphi)$, we find the existence of a global
holomorphic $1$-form $\Phi \in H^0(X\times\cM_H,
\Omega^1_{\cM_H}(End(\cE)))$. Hence, the function
on $\cM_H$ defined by
$$
s \,\mapsto\, \chi(s) \,= \, \int_{X\times\{s\}} g^\ba \tr
(\varphi_\alpha\varphi^*_\ol\beta) g\, dV
$$
is a function of class $\cinf$ on $\cM_H$.

In a similar way, the curvature form $\Omega$ of the
Hermitian-Yang-Mills
connections is a well-defined $End(\cE)$-valued $(1,1)$-form over
$X\times \cM_H$.

For the results of this section, the base space $S$ can be a complex
space (even non-reduced, if necessary). However, in order to
simplify the exposition, we assume smoothness.

Given the projection $X\times S \to S$, where $S$ is also smooth,
the push-forward of a $(n+1,n+1)$-form $\Psi$
(defined on $X\times S$) is a $(1,1)$-form on $S$
given by a fiber integral
$$
\int_{X\times S/S}\Psi, \text{ which we also write as }
\int_{X\times \{s\}} \Psi \text{ or simply } \int_{X}\Psi.
$$

\begin{proposition}\label{prop.fif}
Let $\Omega$ be the curvature form of $(\cE, h)$. Then the following
fiber integral formula holds:
\begin{gather}\label{pwfib}
\omega_{PW} = \frac{1}{2}
\int_X\tr(\Omega\wedge\Omega)\wedge\frac{\omega_X^{n-1}}{(n-1)!}
\hspace{3cm}
\\ \nonumber \hspace{2cm} +
\lambda\int_X\tr(\ii\Omega)\wedge \frac{\omega_X^n}{n!} +\ii
\partial\ol\partial \frac{1}{2}\int_X \tr(\varphi\wedge\varphi^*)
\wedge\frac{\omega_X^{n-1}}{(n-1)!}.
\end{gather}
Here $\lambda$ is determined by
$$
\int_X\tr(\ii\Omega)\wedge\frac{\omega_X^{n-1}}{(n-1)!} =
\lambda\int_X\frac{\omega_X^n}{n!}
$$
is independently of $s\in S$ over any connected component of $S$.
\end{proposition}

Before we prove the proposition, we recall some standard facts.
Concerning Chern character forms, we will use the description as
$$
\ch(\cE,h)= \sum_{k=0}^n \frac{1}{k!}\tr\left( \vtop{\hbox
{$\underbrace{\frac{\ii}{2\pi}\Omega \wedge\dots \wedge
\frac{\ii}{2\pi}\Omega }$}\hbox{\hspace{19mm}{$k$}}} \right)
$$
with
$$
\ch_2(\cE,h)= \frac{1}{2}\left(c_1^2(\cE,h)-2c_2(\cE,h)\right).
$$
In terms of Chern character forms and Chern forms formula
(\ref{pwfib}) reads
\begin{eqnarray}
\frac{1}{4\pi^2} \omega_{PW} &=& - \int_X
\ch_2(\cE,h)\wedge\frac{\omega_X^{n-1}}{(n-1)!} \\ \nonumber& & +
\frac{\lambda}{2\pi} \int_X c_1(\cE,h)\wedge\frac{\omega_X^n}{n!}\\
\nonumber & &+ \frac{\ii}{8\pi^2}\partial\ol\partial\int_X
\tr(\varphi\wedge\varphi^*)\wedge \frac{\omega_X^{n-1}}{(n-1)!}\, .
\end{eqnarray}
Now we will prove the proposition.

\begin{proof}
By definition
$$
\omega_{PW}=\left( \int \tr(R_{i\ol\beta} R_{\alpha\ol\jmath})
g^\ba g\,dV + \int
\tr(\varphi_{\alpha;i}\varphi^*_{\ol\beta;\ol\jmath})\,g\,dV\right)
\ii ds^i\wedge
ds^{\ol\jmath}.
$$
Now
\begin{gather*}
\frac{1}{2}\int_X \tr(\Omega\wedge\Omega)\wedge
\frac{\omega_X^{n-1}}{(n-1)!}  =  - \frac{1}{2}\int_X
\tr(\ii\Omega\wedge\ii\Omega)\wedge \frac{\omega_X^{n-1}}{(n-1)!}\\
\hspace{2cm} = \int\tr(R_{\alpha \ol\jmath}\cdot R_{i\ol\beta} -
R_\ab \cdot R_{i\ol\jmath})g^\ba g \,dV\ii ds^i\wedge ds^\ol\jmath,
\end{gather*}
and from (\ref{hermein}) we have
$$
-\tr(g^\ba R_\ab  \cdot R_{i\ol\jmath})=
\tr(g^\ba[\varphi_\alpha,\varphi^*_\ol\beta]\cdot R_{i\ol\jmath}) -
\lambda \tr R_{i\ol\jmath}.
$$
On the other hand
$$
g^\ba (\varphi_{\alpha;i} \varphi^*_{\ol\beta;\ol\jmath})=  g^\ba
(\varphi_\alpha\varphi^*_\ol\beta )_{i\ol\jmath}
-(\varphi_{\alpha;i\ol\jmath}\cdot \varphi^*_\ol\beta),
$$
and
$$
\varphi_{\alpha;i\ol\jmath}= - \tr [R_{i\ol\jmath},\varphi_\alpha],
$$
from which the claim follows.
\end{proof}

{}From now on, we assume that $X$ is a \ka\ manifold whose
\ka\ form is the Chern form
$$
\omega_X=c_1(\cL,h_\cL)
$$
of a positive hermitian line bundle $(\cL,h_\cL)$. Note that this
implies that $X$ is a complex projective manifold.

Given a proper, smooth holomorphic map $f:\mathcal X \to S$ and a
locally free sheaf $\mathcal F$ on $\mathcal X$, the determinant
line bundle of $\mathcal F$ on $S$ is by definition $\det
\underline{\underline R}f_* \mathcal F$.

The generalized Riemann-Roch theorem by Bismut, Gillet and Soul\'e
\cite{bgs} applies to hermitian vector bundles $(\mathcal F, h)$ on
$\mathcal X$. It states that the determinant line bundle of
$\mathcal F$ on $S$ carries a Quillen metric, whose Chern for equals
the fiber integral
$$
\int_{\mathcal X/S} ch(\mathcal F,h)td(\mathcal X/S,\omega_\mathcal
 X),
$$
where $ch$ and $td$ denote respectively the Chern character form and
the Todd form. (For smooth, proper holomorphic maps over singular
base spaces cf.\ \cite[App.]{f-s}).

We first mention
\begin{equation}
\ch(End(\cE))= r^2 + 2 r \ch_2(\cE) -c_1^2(\cE) + \dots
\end{equation}
so that for the virtual bundle $End(\cE)- \cO^{r^2} $
$$
\ch(End(\cE)- \cO^{r^2}) = 2r ch_2(\cE)-c_1^2(\cE)+ \dots
$$
holds. We use these formulas for hermitian bundles now.
\begin{eqnarray*}
&&\hspace{-8mm} \ch\left(((End(\cE),h)-\cO^{r^2})\otimes \left(
(\cL,h_\cL)-
(\cL^{-1},h_\cL^{-1})\right)^{\otimes (n-1)}\right)  \\
&=& \ch_2\left((End(\cE),h)-\cO^{r^2}\right) \cdot 2^{n-1}
\omega_X^{n-1} +
\ldots\\
&=& \left(2r\left( \frac{1}{2}\tr\left(\frac{\ii}{2\pi}\Omega \wedge
\frac{\ii}{2\pi}\Omega\right)\right) - \left(\tr
\frac{\ii}{2\pi}\Omega    \right)^2      \right)2^{n-1}
\omega_X^{n-1} +\dots\\ &=&
2^{n-1}\left(r\tr\left(\frac{\ii}{2\pi}\Omega\wedge
\frac{\ii}{2\pi}\Omega\right)-\left(\tr
\frac{\ii}{2\pi}\Omega\right)^2\right) \omega_X^{n-1} +\dots
\end{eqnarray*}
The highest exterior power $\Lambda^r \cE$ carries the induced
hermitian metric $\wh h$, for which the following identity holds:
\begin{eqnarray*}
&& \hspace{-2cm}\ch\left(\left(\left(\Lambda^r \cE,\wh
h\right)-(\Lambda^r \cE,\wh h)^{-1}\right)^{\otimes 2} \cdot \left(
(\cL,h_\cL)-(\cL^{-1},h^{-1})\right)^{\otimes(n-1)}\right)\\
\qquad&=& 2^{n+1} c_1^2(\cE,h)\cdot c_1(\cL,h_\cL)^{n-1}  +\dots \\
&=& 2^{n+1} c_1^2(\cE,h)\cdot \omega_X^{n-1} +\dots \\
&=& 2^{n+1}\left(\tr\frac{\ii}{2\pi}\Omega\right)^2 \omega_X^{n-1} +
\dots
\end{eqnarray*}

Hence we have the following theorem:

\begin{theorem}\label{th.fif} The generalized \pw\ form can be
expressed in terms of Chern character forms of hermitian bundles:
\begin{eqnarray*}
&& \hspace{-1cm}\frac{1}{4\pi^2}\omega_{PW} = \\
 && -\frac{1}{2^n r (n-1)!} \int \ch\left((End(\cE)-\cO^{r^2})\otimes
(\cL-\cL^{-1})^{\otimes(n-1)}\right)\\
 &&  - \frac{1}{2^{n+2} r (n-1)!} \int \ch\left(
(\Lambda^r\cE-(\Lambda^r\cE)^{-1})^{\otimes 2}\otimes
(\cL-\cL^{-1})^{\otimes(n-1)}\right)\\
 && + \frac{\lambda}{2\pi}
\frac{1}{2^{n+1} n!}\int \ch\left(\left(
\Lambda^r\cE-(\Lambda^r\cE)^{-1}\right)\otimes
(\cL-\cL^{-1})^{\otimes n}\right)\\
 && + \frac{1}{8\pi^2} \ii
\partial\ol\partial \int \tr(\varphi\wedge \varphi^*)\wedge
\frac{\omega_X^{n-1}}{(n-1)!}\, .
\end{eqnarray*}
\end{theorem}

Let $q: X\times \cM_H \to \cM_H$ be the canonical projection. We
introduce the following determinant line bundles $\delta_j$,
equipped with Quillen metrics $h^Q_j$.
\begin{eqnarray*}
\delta_1 &=& \det \underline{\underline R}q_* \left(
(End(\cE)-\cO^{r^2})\otimes
(\cL-\cL^{-1})^{\otimes(n-1)}\right)\\
 \delta_2 &=& \det
\underline{\underline R}q_* \left(\left(
\Lambda^r\cE-(\Lambda^r\cE)^{-1}\right)^{\otimes 2} \otimes
(\cL-\cL^{-1})^{\otimes(n-1)}\right) \\
  \delta_3 &=& \det
\underline{\underline R}q_* \left(\left(
\Lambda^r\cE-(\Lambda^r\cE)^{-1}\right) \otimes
(\cL-\cL^{-1})^{\otimes n}\right).
\end{eqnarray*}
Setting
$$
\chi=\int \tr (\varphi\wedge \varphi^*)\wedge
\frac{\omega_X^n}{(n-1)!}
$$
we equip the trivial bundle $\cO_{\cM_H}$ with the hermitian metric
$e^\chi$.

\begin{theorem}\label{thm3}
The generalized \pw\ Kähler form is a linear combination of the
$(1,1)$-forms $c_1(\delta_j,h^Q_j)$, $j=1,2,3$, and
$c_1(\cO_{\cM_H},e^\chi)$.
\end{theorem}

\section{A holomorphic closed $2$-form on a
moduli space of Higgs bundles}\label{sec7}

Let $(\cE,\varphi)$ be a universal family of stable Higgs
bundles on a \ka\ manifold $(X,\omega_X)$ over a complex analytic
space $S$ carrying the unique family of Hermitian-Yang-Mills connections
$\{\theta_s\}_{s\in S}$. Using the previous notation, we introduce a
holomorphic two-form $\pi$ on $S$ with $d\pi =0$.

Let $\left.\frac{\partial}{\partial s^i}\right|_{s=s_0}$ be a
tangent vector, and
\begin{gather}\label{rho}
\rho\left(\left.\frac{\partial}{\partial s^i}\right|_{s=s_0}\right)=
\varphi_{\alpha;i}dz^\alpha + R_{i\ol\beta} dz^\ol\beta \hspace{4cm} \\
\nonumber \in C^{1,0}(X\times S, End(\cE_s))\oplus C^{0,1}(X\times
S, End(\cE_s)).
\end{gather}

\begin{definition}
{\rm A two-form
$$
\pi =\pi_{ik}(s)ds^i\wedge ds^k
$$
on $S$ is given by
$$
\pi_{ik}= \int_{X\times\{s\}} \tr(g^\ba(\varphi_{\alpha;i}\cdot
R_{k\ol\beta} - \varphi_{\alpha;k}\cdot R_{i\ol\beta}))g\, dV\, .
$$}
\end{definition}

\begin{lemma}
The two-form $\pi$ is holomorphic, and furthermore, it is of the
form $\pi = d\nu$ for a certain holomorphic $1$-form $\nu$ on $S$.
The forms $\nu$ and $\pi$ descend to the moduli space of Higgs
bundles as holomorphic forms.
\end{lemma}

\begin{proof}
We define $\nu=\nu_i ds^i$ through
$$
\nu_i= 2 \int_{X\times\{s\}}\tr g^\ba \varphi_\alpha R_{i\ol\beta}\,
g\,dV.
$$
Then $d\nu=\pi$ follows immediately, and
\begin{gather*}
\frac{\partial \nu_i}{\partial s^\ol \ell}= \int \tr g^\ba
\varphi_\alpha R_{i\ol\beta;\ol\ell}\, g\, dV
 =- \int \tr g^\ba \varphi_{\alpha;\ol\beta}R_{i\ol\ell}\,g\, dV=0
\end{gather*}
completing the proof of the lemma.
\end{proof}

Denote by $\cM$ the moduli space of stable vector bundles on $X$. As
any stable vector bundle $E$ defines a Higgs bundle $(E,\varphi)$
with Higgs field $\varphi=0$. The \he\ connection on the stable
vector bundle $E$ coincides with the Hermitian-Yang-Mills connection
on $(E, 0)$. We have an embedding $i: \cM \hookrightarrow \cM_H$
into the corresponding moduli space of stable Higgs bundles. Let
$\cM^s_H \, \subset\, \cM_H$ denote the Zariski open subset defined
by all Higgs bundles $(E,\varphi)$ with $E$ stable. Therefore, we
have a retraction $f: \cM^s_H \to \cM$ that sends any $(E,\varphi)$
to $E$.

\begin{proposition}
The forms $\nu$ and $\pi$ vanish on the fibers of $f$.
\end{proposition}

\begin{proof}
On the level of base spaces of universal deformations, $f$ is a
submersion of the form $pr_1:S=S'\times S'' \to S'$. Let $s=(s',s'')
\in S$. Let $v=\frac{\partial}{\partial s^i}\in T_sS$ be a tangent
vector with
$$
\rho(v)= \varphi_{\alpha;i}dz^\alpha + R_{i\ol\beta}dz^\ol\beta
$$
in the sense of (\ref{rho}). Then $pr_1(v)$ is represented by
$R_{i\ol\beta} dz^\ol\beta$, and if $v\in T_sS'' \hookrightarrow T_s
S$, the form $R_{i\ol\beta} dz^\ol\beta$ is $\ol\partial$-exact.
Conversely let $R_{k\ol\beta}=E_{;\ol\beta}$ for some section $E$ of
$End(\cE|X\times S'')$. Then
$$
\nu_k(s)= \int \tr g^\ba \varphi_\alpha E_{;\ol\beta} g\, dV=0
$$
proving the proposition.
\end{proof}

If the form $\nu$ is nonzero on $\cM^s_H$, then it does not
descend under $f$.

\section{Non-Abelian Hodge symmetry, symplectic and hyper-Kähler
structure}\label{nonabhodge}

Let $(E,\varphi)$ be a stable Higgs bundle over $X$ equipped with a
Hermitian-Yang-Mills connection. We first provide the space
of infinitesimal deformations of $(E,\varphi)$, namely
${\mathbb H}^1(C^{\bullet\bullet})$,
with a quaternionic structure under an assumption on
$\mathbb H^2(C^{\bullet\bullet})$.

\noindent {\bf Assumption.} {\it For the rest of this section we
restrict ourselves to stable Higgs bundles $(E,\varphi)$ satisfying
the following two conditions:
\begin{description}
\item[A] The rational characteristic classes of the projective bundle
${\mathbb P}(E)$ over $X$ are assumed to be zero. This is equivalent
to the condition that the Hermitian-Yang-Mills connection on
$(E,\varphi)$ is projectively flat, i.e.,
\begin{equation}
\label{projflat}R_\ab + [\varphi_\alpha,\varphi^*_\ol\beta]=
\lambda\cdot g_\ab \cdot {\rm id}_E
\end{equation}
for some $\lambda \in {\mathbb C}$.
\item[B] $\dim \mathbb H^2(C^{\bullet\bullet})\,=\, 1$.

In view of Lemma \ref{lale} this is equivalent to the condition that
\begin{equation}\label{h2}
\mathbb H^2(C^{\bullet\bullet})\,=\, \mathbb C\cdot \omega_X\cdot
{\rm id}_E\, .
\end{equation}
\end{description}
}

The above condition $\bf A$ can be replaced by
\begin{description}
\item[\bf A'] The Higgs field $\varphi$ is
$\partial_\theta$-closed with
respect to the Hermitian-Yang-Mills connection, i.e.\
\begin{equation}\label{aprime}
\varphi_{\alpha;\gamma}=\varphi_{\gamma;\alpha}.
\end{equation}
\end{description}

\noindent{\bf Involution.} We have an involution of the space of
$End(E)$-valued $1$-forms defined by
\begin{eqnarray}\label{iota}
\iota: C^1 &\to& C^1\\
(a,b)&\mapsto& (-b^*, a^*)\nonumber
\end{eqnarray}
where $C^1$ is defined in \eqref{de-cr}.
Obviously $\iota^2= - {\rm id}_{C^1}$. We shall see this
involution descends to
the space of infinitesimal deformations of $(E,\varphi)$.

\begin{proposition}
Let $\eta \in C^1$ be harmonic (so $\eta$ gives an
infinitesimal deformation of $(E,\varphi)$). Then $\iota(\eta)$ is
harmonic.
\end{proposition}
\begin{proof}
Let $\eta\neq 0$. Since deformations are not obstructed, we can
assume that there is a coordinate system on the base $S$ of a
universal deformation so that $\eta$ is of the form $\eta_i$ in the
sense of (\ref{defetai}), i.e.\
$$
\iota\eta_i=(-R_{\alpha\ol\imath}
dz^\alpha,\varphi^*_{\ol\beta;\ol\imath}dz^\ol\beta)\, .
$$
\medskip
{\it Claim 1.}
$$
d(\iota\eta_i)= \xi_{;\ol\imath}\, ,
$$
where
$$
\xi=(\varphi_{\alpha;\gamma}dz^\alpha\wedge dz^\gamma,
(R_\ab+[\varphi_\alpha,\varphi^*_\ol\beta])dz^\alpha\wedge
dz^\ol\beta, -\varphi^*_{\ol\beta;\ol\delta}dz^\ol\beta\wedge
dz^\ol\delta)\, .
$$
{\it Proof of Claim 1.} Let
$d(\iota\eta_i)=(d(\iota\eta_i)_1,d(\iota\eta_i)_2,d(\iota\eta_i)_3)$.
Then
$$
d(\iota\eta_i)_1 = [R_{\alpha\ol\imath}dz^\alpha, \varphi_\gamma
dz^\gamma] = \varphi_{\alpha;\gamma\ol\imath} dz^\alpha\wedge
dz^\gamma\, .
$$
Next,
\begin{eqnarray*}
d(\iota\eta_i)_2& =& \db(-R_{\alpha\ol\imath}dz^\alpha) +
[\varphi^*_{\ol\beta;\ol\imath}dz^\ol\beta,\varphi_\alpha
dz^\alpha]\\
&=& (R_\ab +[\varphi_\alpha,\varphi^*_{\ol\beta;
\ol\imath}])_{;\ol\imath}\;dz^\alpha\wedge
dz^\ol\beta\, ,
\end{eqnarray*}
and
$$
d(\iota\eta_i)_3= \db(\varphi^*_{\ol\beta;\ol\imath}dz^\ol\beta) = -
\varphi_{\ol\beta;\ol\delta\ol\imath} dz^\ol\beta \wedge
dz^\ol\delta
$$
proving Claim 1.

{\it Claim 2.}
$$
d^*(\iota \eta_i)=0\, .
$$
{\it Proof of Claim 2.}
\begin{eqnarray*}
d^*(\iota\eta_i)&=&
+\Lambda[R_{\alpha\ol\imath}dz^\alpha,\varphi^*_\ol\beta
dz^\ol\beta]+\db^*(\varphi^*_{\ol\beta;\ol\imath} dz^\ol\beta)\\ &=&
g^\ba([R_{\alpha\ol\imath},\varphi^*_\ol\beta ]-
\varphi^*_{\ol\beta;\ol\imath \alpha} )=0\, .
\end{eqnarray*}

{\it Claim 3.}
$$
d\xi=0.
$$
{\it Proof of Claim 3.} Let
$$
d\xi=(d\xi_1, \ldots,d\xi_4)
$$
where $d\xi_i$ is the $i$-th component of $d\xi$
in the decomposition in \eqref{de-cr}. Then
$$
d\xi_1= [-\varphi_{\alpha;\gamma}dz^\alpha \wedge dz^\gamma,
\varphi_\sigma dz^\sigma]=0
$$
because of (\ref{integrability}).

Next,
\begin{eqnarray*}
d\xi_2&=&\db(\varphi_{\alpha;\gamma}dz^\alpha\wedge dz^\gamma) +
[(R_\ab +[\varphi_\alpha,\varphi^*_\ol\beta])dz^\alpha\wedge
dz^\ol\beta,\varphi_\gamma
dz^\gamma]\\
&=& (\varphi_{\alpha;\gamma\ol\beta}
-[R_\ab,\varphi_\gamma]-[[\varphi_\alpha,\varphi^*_\ol\beta],\varphi_\gamma])
dz^\alpha\wedge dz^\gamma \wedge dz^\ol\beta\\
&=& ([\varphi_\alpha,R_{\gamma\ol\beta}]
+[\varphi_\gamma,R_\ab]-[[\varphi_\alpha,\varphi^*_\ol\beta],
\varphi_\gamma])dz^\alpha\wedge
dz^\gamma\wedge dz^\ol\beta.
\end{eqnarray*}
The first two terms together are symmetric in $\alpha$ and $\gamma$,
also the third term because of (\ref{integrability}) so that
$d\xi_2$ vanishes.

Finally,
$$
d\xi_3=
-[-\varphi^*_{\ol\beta;\ol\delta},\varphi_\alpha]dz^\ol\beta\wedge
dz^\ol\delta\wedge dz^\alpha + \db(R_\ab
+[\varphi_\alpha,\varphi^*_\ol\beta])dz^\alpha \wedge dz^\ol\beta=0,
$$
and
$$
d\xi_4=\db(-\varphi^*_{\ol\beta;\ol\delta}dz^\ol\beta\wedge
dz^\ol\delta)=0
$$
proving Claim 3.

{\it Claim 4.}
$$
d^*\xi=0.
$$
{\it Proof of Claim 4.} Let $\xi=(u,v,w)$. Here projective flatness will
be used: According to (\ref{dstar1}) the first component of $d^*\xi$
equals
$$
d^*\xi_1= \Lambda([\varphi_{\alpha;\gamma}dz^\alpha\wedge dz^\gamma,
\varphi^*_\ol\delta dz^\ol\delta]) + \dbs v\, ,
$$
where $d^*\xi_i$ is the $i$-th component of $d^*\xi$
in the decomposition in \eqref{de-cr}.
Because of \eqref{projflat} we have $\dbs v=0$, and
$$
d^*\xi_1
=\frac{1}{2}g^{\ol\beta\gamma}[(\varphi_{\alpha;\gamma}-
\varphi_{\gamma;\alpha}),\varphi^*_\ol\beta]dz^\alpha=
\frac{1}{2}g^{\ol\beta\gamma} (R_{\ab;\gamma} -
R_{\gamma\ol\beta;\alpha} )dz^\alpha =0\, .
$$
We compute $d^*\xi_2$. By \eqref{projflat}, the term
$\wt\Lambda[v,\varphi^*]$ vanishes, and
\begin{gather*}
\dbs(-\varphi^*_{\ol\beta;\ol\delta} dz^\ol\beta\wedge dz^\ol\delta)
=
\frac{1}{2}g^{\ol\delta\gamma}(\varphi^*_{\ol\beta;\ol\delta\gamma}
- \varphi^*_{\ol\delta;\ol\beta\gamma}) dz^\ol\beta =\\
\frac{1}{2}g^{\ol\delta\gamma}(-[\varphi^*_\ol\beta,R_{\gamma\ol\delta}]+
[\varphi^*_\ol\delta,R_{\gamma\ol\beta}])dz^\ol\beta=
\frac{1}{2}g^{\ol\delta\gamma}([[\varphi^*_\ol\beta,\varphi_\gamma],
\varphi^*_\ol\delta
]-[[\varphi^*_\ol\delta,\varphi_\gamma],\varphi^*_\ol\beta ]
)dz^\ol\beta\, .
\end{gather*}
This form vanishes because of the integrability condition
\eqref{integrability}, proving Claim 4. The computation is the same
with {\bf A} replaced by {\bf A'}.

Now, by our above assumption {\bf B}, the harmonic form $\xi$ is of
the form $c(s)\cdot \omega_X {\rm id}_E$. From the definition of $\xi$
and \eqref{projflat} we know that $c(s)=\lambda$ is independent of $s\in
S$. Hence $d(\iota \eta_i)=0$. This completes the proof of the
proposition.
\end{proof}

\begin{corollary}
The $d$-closed holomorphic $2$-form $\pi$ is non-degenerate at. In
particular, the dimension of $\mathbb H^1(C^{\bullet\bullet})$ is
even.
\end{corollary}

\begin{proof}
In the above notation,
$$
\pi_{ik}=\langle \eta_k, \iota( \eta_i) \rangle.
$$
As $\iota$ takes harmonic elements of $C^1$ to harmonic elements of
$C^1$, it induces a bijective map of $\mathbb H^1$ to itself. This
shows the non-degeneracy.
\end{proof}

\begin{remark}
For $X = {\mathbb C}{\mathbb P}^2$, there are examples of
stable Higgs bundles where
$\dim {\mathbb H}^1(C^{\bullet\bullet})$ is odd. (They do not
satisfy the two assumptions.)
\end{remark}

We insert an observation about quaternions: Assume that $(\mathbb H,
I,J,K)$ is identified with $\C^2=\{(z,w)\}$ in such a way that
$I(z,w)=(i z,i w)$, $i=\ii$. Namely, we identify $(z,w)=(x+i y,u+i
v)$ with $x+Iy + (u+I v)J$. Then $J(z,w)=(-\ol w,\ol z)$, and
$K(z,w)= IJ(z,w)=(-i \ol w, i \ol z)$.

This suggests
the following almost quaternionic structure on
${\mathbb H}^1(C^{\bullet\bullet})$: For
$\eta=(a,b)\in {\mathbb H}^1(C^{\bullet\bullet})$
we have $\iota (a,b)=(-b^*,a^*)$, and
hence for a tangent vector
$$
I(\eta)= \sqrt{-1}\cdot \eta,\qquad  J(\eta)=\iota (\eta),
\qquad K(\eta)= \sqrt{-1} \cdot \iota(\eta).
$$
By definition the equations for an almost quaternionic structure are
verified.

For complex tangent vectors $\eta,\vartheta$ of the base space have
\begin{eqnarray*}
\langle I\eta,\vartheta\rangle&=&\omega^{PW}(\eta,\vartheta),\\
\langle J\eta,\vartheta\rangle &=& \ol\pi(\eta,\vartheta),\\
\langle K\eta,\vartheta\rangle &=& i\ol\pi(\eta, \vartheta),
\end{eqnarray*}
where $\ol\pi$ denotes the conjugate of $\pi$. The corresponding
differential forms are closed (and non-degenerate), so that together
with \cite[Lemma (6.8)]{hitchin} the following holds:

\begin{theorem}\label{tf}
Consider the Zariski open subset $W$ of the moduli of stable Higgs
bundles over which the condition $\mathbb H^2 = \C\cdot \omega_X
{\rm id}_E$ holds. Assume that either the Hermitian-Yang-Mills
connections are projectively flat or that the Higgs fields are
$\partial_\theta$-closed with respect to the Hermitian-Yang-Mills
connection. Then $W$ carries a natural hyper-Kähler structure,
related to the \pw structure $\omega^{PW}$ and the holomorphic
symplectic form $\pi$.
\end{theorem}

We note that our proof shows more.
\begin{proposition}
Consider the moduli space of Higgs bundles of the form
$(E,\varphi)$, where the determinant bundle of $E$ is a fixed line
bundle, and where $\rm{trace} (\varphi)=0$. This moduli space
carries a natural hyper-Kähler structure, related to the \pw
structure $\omega^{PW}$ and the holomorphic symplectic form $\pi$,
provided condition {\bf A} or condition {\bf A'} holds along with
the following weaker form of condition {\bf B}:
\begin{description}
  \item[B'] $$\mathbb H^2(D^\bullet_0)=0,$$
  where $D^\bullet_0$ is the complex in Lemma~\ref{decomp}.
\end{description}
\end{proposition}

\noindent
{\bf Acknowledgements.} The first named author
wishes to thank the University of Marburg for its hospitality
and support, while
the second named author would like to thank the Tata Institute for
Fundamental Research for hospitality and support. This work was
partially supported by DFG-Schwerpunktprogramm ''Globale Methoden in
der Komplexen Geometrie''.



\begin{thebibliography}{O-T-T2}

\bibitem[B-G-S]{bgs}
Bismut, J.-M.; Gillet, H.; Soul\'e, C.: Analytic torsion and
holomorphic determinant bundles. I: Bott-Chern forms and analytic
torsion. II: Direct images and Bott-Chern forms. III: Quillen
metrics on holomorphic determinants. Commun. Math. Phys.
\textbf{115}, 49--78, 79--126, 301--351 (1988)

\bibitem[Bi]{biswas} Biswas, I.:
A remark on a deformation theory of Green and Lazarsfeld. Jour.
Reine Angew. Math. \textbf{449}, 103--124 (1994)

\bibitem[B-R]{br} Biswas, I.; Ramanan, S.: An infinitesimal study
of the moduli of Hitchin pairs. Jour. Lond. Math. Soc. \textbf{49},
219--231 (1994)

\bibitem[B-S]{b-s2} Biswas, I.; Schumacher, G.: Determinant bundle,
Quillen metric, and Petersson-Weil form on moduli spaces. Geom.
Funct. Anal. \textbf{9}, 226--256 (1999)

\bibitem[F-S]{f-s} Fujiki, A., Schumacher, G.: The moduli space of compact
extremal Kähler manifolds and generalized Petersson-Weil metrics.
Publ.\ RIMS, Kyoto Univ.\ {\bf 26} 101--183 (1990).

\bibitem[Hi]{hitchin} Hitchin, N.: The self-duality equations on a
Riemann surface. Proc. Lond. Math. Soc. \textbf{55}, 59--126 (1987)

\bibitem[Sch]{sch}
Schumacher, G.: Moduli as algebraic spaces. Complex analysis in
several variables---Memorial Conference of Kiyoshi Oka's
Centennial Birthday, 283--288, Adv. Stud. Pure Math., 42, Math.
Soc. Japan, Tokyo, 2004.

\bibitem[S-T]{s-t}Schumacher, G.; Toma, M.: On the
Petersson-Weil metric for the moduli space of Hermite-Einstein
bundles and its curvature. Math. Ann. \textbf{293}, 101--107 (1992).

\bibitem[Si1]{simpson0} Simpson, C.T.: Constructing variations
of Hodge structure using Yang-Mills theory and applications
to uniformization. Jour. Amer. Math. Soc. \textbf{1}, 867--918 (1988)

\bibitem[Si2]{simpson} Simpson, C.T.: Higgs bundles and local
systems. Publ. Math. I.H.\'E.S. \textbf{75}, 5--95 (1992)

\bibitem[Si3]{simpson2} Simpson, C.T.: Moduli of
representations of the fundamental group of a smooth
projective variety. II. Publ. Math. I.H.\'E.S.
\textbf{80}, 5--79 (1994)

\bibitem[Va]{va} Varouchas, J.: Stabilit\'e de la classe
des vari\'et\'es K{ä}hl\'eriennes par certains morphismes propres.
Invent.\ Math. \textbf{77} (1984), 117--127.

\end{thebibliography}
\end{document}